\pdfoutput=1
\RequirePackage{ifpdf}
\ifpdf 
\documentclass[pdftex]{sigma}
\else
\documentclass{sigma}
\fi

\begin{document}

\allowdisplaybreaks

\renewcommand{\PaperNumber}{029}

\FirstPageHeading

\ShortArticleName{Orbit Representations from Linear mod 1 Transformations}

\ArticleName{Orbit Representations\\ from Linear mod 1 Transformations}

\Author{Carlos CORREIA RAMOS~$^\dag$, Nuno MARTINS~$^\ddag$ and Paulo R.~PINTO~$^\ddag$}

\AuthorNameForHeading{C.~Correia Ramos, N.~Martins and P.R.~Pinto}

\Address{$^\dag$~Centro de Investiga\c c\~ao em Matem\'atica e Aplica\c
c\~oes,\\
\hphantom{$^\dag$}~R.\ Rom\~{a}o Ramalho, 59, 7000-671 \'{E}vora, Portugal}
\EmailD{\href{mailto:ccr@uevora.pt}{ccr@uevora.pt}}

\Address{$^\ddag$~Department of Mathematics, CAMGSD,
Instituto Superior T\'{e}cnico,\\
\hphantom{$^\ddag$}~Technical University of Lisbon,
Av.\  Rovisco Pais, 1049-001 Lisboa, Portugal}
\EmailD{\href{mailto:nmartins@math.ist.utl.pt}{nmartins@math.ist.utl.pt}, \href{mailto:ppinto@math.ist.utl.pt}{ppinto@math.ist.utl.pt}}

\ArticleDates{Received March 14, 2012, in f\/inal form May 09, 2012; Published online May 16, 2012}

\Abstract{We show that every point $x_0\in [0,1]$ carries a representation
of a $C^*$-algebra that encodes the orbit structure of the
linear mod 1 interval map $f_{\beta,\alpha}(x)=\beta x +\alpha$. Such $C^*$-algebra is generated
by partial isometries arising from the subintervals of monotonicity of the underlying map $f_{\beta,\alpha}$.
 Then we prove that such representation is irreducible.
Moreover two such of representations are unitarily equivalent
if and only if the points belong to the same generalized orbit,
for every $\alpha\in [0,1[$ and $\beta\geq 1$.}

\Keywords{interval maps; symbolic dynamics; $C^*$-algebras; representations of algebras}

\Classification{46L55; 37B10; 46L05}

\section{Introduction}

A famous class of representations of the Cuntz algebra
$\mathcal{O}_n$ called permutative representations were studied
and classif\/ied by Bratteli and Jorgensen in \cite{BJ,BJ2}. From
the applications viewpoint, and besides its own right,
applications of representation theory of Cuntz and Cuntz--Krieger
algebras to wavelets, fractals, dynamical systems, see e.g.\
\cite{BJ,BJ2,marcolli}, and quantum f\/ield theory in~\cite{AK} are
particularly remarkable. For example, it is known that these
representations of the Cuntz algebra serve as a computational tool
for wavelets analysts, see~\cite{Palle2005}. This is clear because
such a representation on a Hilbert space $H$ induces a subdivision
of $H$ into orthogonal subspaces. Then the problem in wavelet
theory is to build orthonormal bases in $L^2(\mathbb{R})$ from
these data. Indeed this can be done \cite{Dau} and these wavelet
bases have the advantages over the earlier known basis
constructions (one advantage is the ef\/f\/iciency of computation).
This method has also been applied to the context of fractals that
arise from af\/f\/ine iterated function systems \cite{DuPalle}. Some
of these results have been extended to the more general class of
Cuntz--Krieger algebras, see \cite{CMP, marcolli} and {\it subshift
$C^*$-algebras} \cite{toke2, exel,matsumoto1} (whose underlying
subshift is not necessarily of f\/inite type) in~\cite{CMP6}.

Symbolic dynamics is one of the main tools that we have used in
\cite{CMP6, CMP} to construct representations of Cuntz,
Cuntz--Krieger and subshift $C^*$-algebras. The $C^*$-algebra is naturally
associated to the given interval map and the Hilbert spaces
naturally arise from the generalized orbits of the interval map.
For a particular family of interval maps, we were able to recover
Bratteli and Jorgensen permutative representations in \cite{CMP}
among the class of Markov maps (which underline the Cuntz--Krieger
algebras of the transition matrix). We remark that while these
Cuntz and Cuntz--Krieger algebras are naturally associated to the
so-called Markov or periodic dynamical systems, the subshift
$C^*$-algebras are ready to incorporate bigger classes of interval
maps.

The interval maps that we treat in \cite{CMP6} are unimodal maps
(that have precisely two subintervals of monotonicity). Then the
representations of the subshift $C^*$-algebra constructed in
\cite{CMP6} are shown to coincide with the ones constructed in
\cite{CMP} from the Cuntz--Krieger algebra, provided the underlying
dynamical system is periodic and therefore has a f\/inite transition
Markov Matrix. However, the proof of the irreducibility of the
subshift $C^*$-algebras representations (for unimodal maps without a
f\/inite transition Markov matrix) rely on the structure of these
unimodal interval maps where the $C^*$-algebra is generated by two
partial isometries.

In this paper we construct representations (of a subshift
$C^*$-algebra generated by $n$ partial isometries) from a family of
interval maps and prove the irreducibility of these
representations (avoiding the unimodal maps techniques used in
\cite{CMP6}).

\begin{figure}[htp]
\centering
  \includegraphics[width=2in]{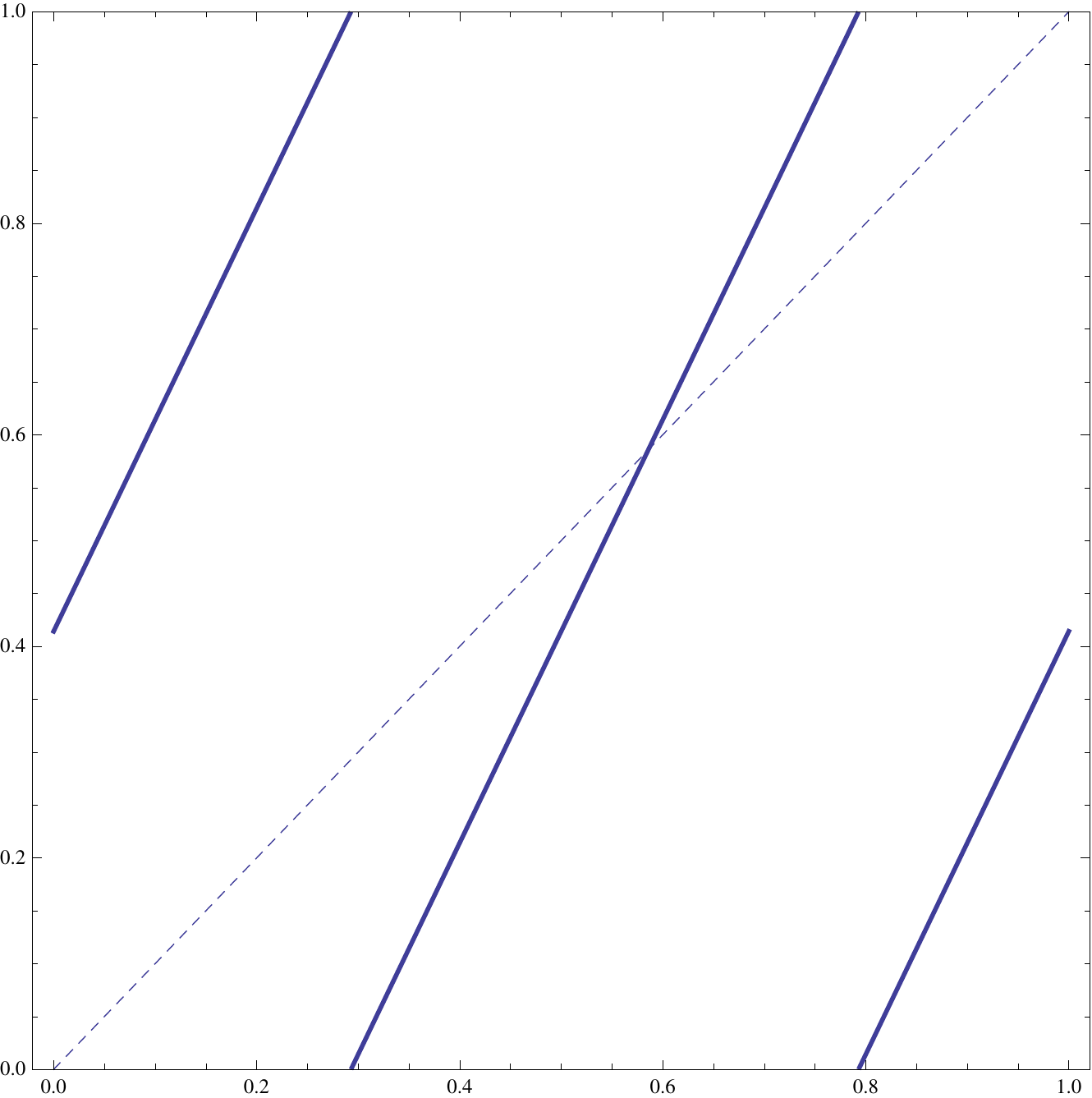}\\
  \caption{Graph of $f_{\beta, \alpha}$ with $ \alpha =
 \sqrt{2}-1$ and $\beta=2\ (n=3)$.} \label{fig2}
\end{figure}

Namely, we yield and study representations of a certain $C^*$-algebra
on the generalized orbit $\bigcup_{i\in\mathbb{Z}}
f_{\beta,\alpha}^j(x_0)$ of every point $x_0\in [0,1]$ from the
interval map $f_{\beta,\alpha}: [0,1]\to[0,1]$ def\/ined by
\begin{equation}
f_{\beta,\alpha}(x)=\beta x+\alpha \quad (\mathrm{mod}\ 1)\quad
\hbox{with}\quad \beta \geq 1\quad \hbox{and}\quad \alpha \in
\lbrack 0,1[, \label{linearmod1}
\end{equation}
by f\/ixing the parameters $\alpha$ and $\beta$. The underlying
$C^*$-algebra $\mathcal{O}_{\Lambda_{f_{\beta,\alpha}}}$
is generated by $n$ partial isometries where $n$ is
the number of monotonicity subintervals of $f_{\beta,\alpha}$.
See Fig.~\ref{fig2} for a graph of one such map.

We show that the representation is irreducible.
Moreover the representations of the same algebra on the orbits
of the two points $x_0$ and $y_0$ are unitarily equivalent if
and only if the orbits coincide.
If the parameters $\alpha$ and $\beta$ are so that the dynamical system
$([0,1], f_{\beta,\alpha})$ is periodic,
the above results were obtained by  \cite{CMP3, CMP},
where the relevant $C^*$-algebra is the
Cuntz--Krieger~$\mathcal{O}_{A_{f_{\beta,\alpha}}}$
and~$A_{f_{\beta,\alpha}}$ is the underlying Markov transition 0-1
(f\/inite) matrix of $f_{\beta,\alpha}$. A~further particular case
is obtained when $\beta=n$ is an integer and $\alpha=0$, in which
case the $n\times n$ matrix $A_{f_{n,0}}=(a_{i,j})$ is full,
$a_{ij}=1$ for all $i$, $j$, and thus recovering the Cuntz algebra
$\mathcal{O}_n$ representations yielded in~\cite{BJ,BJ2} using
wavelet theory framework.

Of course we may fairly easily prove that we do get representations of such
$C^*$-algebras in the context of piecewise monotone maps, but the main concern is
how to show irreducibility (and unitarily equivalence) of such the representations,
thus giving a rich family of representations attached to interval maps.
The periodic or non-periodic cases for which we have  $n=2$ subintervals of monotonicity
were carried out in~\cite{CMP6}. We generalize here the construction of the representations
for generic piecewise monotone interval maps and prove irreducibility and unitarily
equivalence for the dynamical systems arising from equation~(\ref{linearmod1}) above~--  and obviously
we look for values of the parameters~$\alpha$ and~$\beta$ for which we do not have a f\/inite Markov
transition matrix for~$f_{\beta,\alpha}$ and thus we really get representations that cannot be
recovered from~\cite{CMP6,CMP3, CMP}.

A more detailed description of the paper is as follows. In Section~\ref{prel} we provide some background material f\/irst on the
operator algebras setup and then in symbolic dynamics~\cite{M-T}.
The main results are in Section~\ref{matexel}. We consider a partition $\mathcal{I}$ of the
interval $I=[0,1]$ into subintervals so that the restriction of
$f_{\beta,\alpha}$ to each of these subintervals is monotone. Then
for every $x_0\in I$ we explicitly def\/ine in equation~(\ref{matsumotot}) a
linear operator on the Hilbert space~$H_{x_0}$ that arises from
the generalized orbit of~$x_0$, for every such subinterval. The Hilbert
space~$H_{x_0}$ encodes the generalized orbit orbit$(x_0)$ of $x_0$ and in
fact every  $\xi\in \hbox{orbit}(x_0)$ is regarded as a vector $|\xi\rangle$ in
$H_{x_0}$ using Dirac's notation.
Then we prove that these linear operators do satisfy the relation they
ought to satisfy, leading to a representation~$\rho_{x_0}:
\mathcal{O}_{A_{f_{\beta,\alpha}}}\to B(H_{x_0})$ of the
$C^*$-algebra $\mathcal{O}_{A_{f_{\beta,\alpha}}}$, as in Proposition~\ref{matrep}.

Then the main result of this paper is Theorem~\ref{nicethm11} where we
show that representation $\rho_{x_0}$ of the underlying $C^*$-algebra
$\mathcal{O}_{A_{f_{\beta,\alpha}}}$ is irreducible and that
two such representations $\rho_{x_0}$ and $\rho_{y_0}$ are
unitarily equivalent if and only if $y_0$ belongs to the
generalized orbit of $x_0$. The new ingredient involved in the
proof of Theorem \ref{nicethm11} is the computation of the commutant $\mathcal{A}_{\beta,\alpha}^\prime$ of $C^*$-algebra~$\mathcal{A}_{\beta,\alpha}$ (generated by the operators def\/ined in equations~(\ref{v1n}) and
(\ref{u1n})) in $B(H_{x_0})$. Indeed, as soon as we prove that the commutant is trivial $\mathcal{A}_{\beta,\alpha}^\prime =\mathbb{C}{\bf 1}$, we only have to show that $\mathcal{A}_{\beta,\alpha}^\prime$ contains $\rho_{x_0}(\mathcal{O}_{A_{f_{\beta,\alpha}}})^\prime$ as in Proposition~\ref{pro23}.

\section{Preliminaries}
\label{prel}

In this section we provide some necessary background, starting with the
operator algebras we obtain from dynamical systems. A
representation of a $\ast$-algebra $\mathcal{A}$ on a complex
Hilbert space~${H}$ is a $\ast$-homomorphism
$\pi:\mathcal{A}\rightarrow{B}(H)$ into the $\ast $-algebra
${B}({H})$ of bounded linear operators on~${H}$. Usually
representations are studied up to unitary equivalence. Two
representations $\pi :\mathcal{A} \rightarrow {B}({H})$ and
$\widetilde{\pi }:\mathcal{A} \rightarrow {B}(\widetilde{{H}})$
are (unitarily) equivalent if there is a unitary operator $U:
H\rightarrow \widetilde{{H}}$ (i.e., $U$ is a surjective isometry)
such that
\[
U\pi (a)=\widetilde{\pi }(a)U \quad \text{for every} \ a\in
\mathcal{A},
\]
and in this case we write $\pi\sim\tilde{\pi}$. A representation
$\pi: \mathcal{A}\rightarrow {B}({H})$ of some $\ast $-algebra is
said to be \textit{irreducible} if there is no non-trivial
subspace of ${H}$ invariant with respect to all operators $\pi
(a)$ with $a\in \mathcal{A}$. A well known result, see e.g.\
\cite[Proposition 3.13]{Ped}, says that $\pi $ is irreducible if
and only if
\begin{equation}
x\in {B}({H}):\ x\pi (a)=\pi (a)x \quad \text{for\ all} \ a\in
\mathcal{A}\ \Longrightarrow \ x=\lambda {\bf 1},  \label{irrr}
\end{equation}
for some complex number $\lambda$, where ${\bf 1}$ denotes the identity
of $B(H)$. By the very def\/inition of commutant, (\ref{irrr})
can be restated as follows: $\pi (\mathcal{A})^{\prime
}=\mathbb{C}{\bf 1}$.
We will be interested in some classes of~$C^\ast$-algebras (=
Banach $*$-algebras such that $||aa^\ast||=||a||^2$ holds for all
$a$, see e.g.~\cite{Ped}). Besides, if we have a representation
$\pi: A \to {B}({H})$ of a~$C^\ast$-algebra $A$, then $\pi$ being
a~$*$-homomorphism implies that $||\pi(a)||\leq ||a||$ for all $a\in
A$, thus $\pi$ is automatically continuous, see
also e.g.~\cite[Secttion~1.5.7]{Ped}).

\subsubsection*{Subshift $\boldsymbol{C^*}$-algebras}

Let $\Lambda\subseteq \Sigma^\mathbb{N}$ be a subshift with a
f\/inite alphabet $\Sigma=\{1,\dots,n\}$. Exel \cite{exel} and
Matsumoto \cite{matsumoto1} constructed $C^*$-algebras associated to
$\Lambda$. Carlsen and Silvestrov \cite{toke2} unif\/ied the two constructions that led them to a
$C^*$-algebra $\mathcal{O}_{\Lambda}$, which is unital and generated
partial isometries $\{t_i\}_{i\in \Sigma}$. Then the partial
isometries that generate $\mathcal{O}_{\Lambda}$ obey the
following relations:
\begin{gather*} \sum t_i t_i^\ast = {\bf 1},\qquad
t_{\alpha }^{\ast}t_{\alpha }t_{\beta } =
t_{\beta }t_{\alpha \beta }^{\ast}t_{\alpha \beta },\qquad
t_\alpha^\ast t_\alpha t_\beta^\ast t_\beta = t_\beta^\ast t_\beta
t_\alpha^\ast t_\alpha,
 \end{gather*}
 where $t_\alpha=t_{\alpha_1}\cdots t_{\alpha_{|\alpha|}}$ and
$t_\beta=t_{\beta_1}\cdots t_{\beta_{|\beta|}}$ with $\alpha$, $\beta$ admissible words (if $\alpha=(\alpha_1,\dots,\alpha_k)$ with $\alpha_i\in \Sigma$ we denote by $|\alpha|$ the length $k$ of $\alpha$).
The algebra $\mathcal{O}_{\Lambda}$ is called
the $C^*$-algebra associated to the subshift $\Lambda$ or subshift
$C^*$-algebra.
Important properties of the subshift $C^\ast$-algebra
$\mathcal{O}_{\Lambda}$ (e.g.\ simplicity) are naturally inherited
from properties of the subshift~$\Lambda$.
 If
$\Lambda$ is a subshift of f\/inite type, then ${\mathcal
O}_\Lambda$ is nothing but the well known Cuntz--Krieger algebra
\cite{matsumoto1}, where {\it the} Cuntz--Krieger algebra
$\mathcal{O}_{A}$ associated to a 0-1 matrix $A=(a_{ij})$ is the
$C^{\ast }$-algebra \cite{CK} generated by (non-zero) partial
isometries $s_{1},\dots,s_{n}$ satisfying:
\begin{equation}
s_{i}^{\ast }s_{i}=\sum_{j=1}^{n}a_{ij}s_{j}s_{j}^{\ast }\quad
(i=1, \dots,n),\qquad \sum_{i=1}^{n}s_{i}s_{i}^{\ast }={\bf 1},
\label{ckrelns}
\end{equation}
and the Cuntz algebra $\mathcal{O}_n$ is the Cuntz--Krieger algebra
$\mathcal{O}_A$ with $A$ full $a_{ij}=1$ for all $1\leq i,j\leq
n$.

\subsection{Symbolic dynamics on piecewise monotone interval maps}

Let $f:I\rightarrow I$ be a piecewise monotone map of the interval
$I$ into itself, that is, there is a~minimal partition of open
sub-intervals of~$I$, $\mathcal{I}=\{I_{1},\dots,I_{n}\}$ such that
$\overline{\bigcup _{j=1}^{n}I_{j}}=I$ and $f_{|I_{j}}$ is continuous
monotone, for every $j=1,\dots,n$. We def\/ine $f_{j}:=f_{|I_{j}}$.
The inverse branches are denoted by $f_{j}^{-1}: f(I_{j})
\rightarrow I_{j}$. Let $\chi _{I_{i}}$ be the characteristic
function on the interval $I_{i}$. The following are naturally
satisf\/ied
\[
f\circ f_{j}^{-1}(x)=\chi _{f(I_{j})}(x)x, \qquad f_{j}^{-1}\circ f_{|I_{j}}(x)=\chi _{I_{j}}(x)x.
\]

Let $\{1,2, \dots,m\}$ be the alphabet associated to some
partition $\{J_{1},\dots,J_{m}\}$ of open sub-intervals
of~$I$  so that $\overline{\bigcup _{j=1}^{m}J_{j}}=I$, not
necessarily~$\mathcal{I}$. The \textit{address map, }is def\/ined by
\[
{\rm ad}: \ \bigcup_{j=1}^{m}J_{j}\rightarrow \{1,2,\dots,m\},\qquad
{\rm ad}(x)=i \quad \text{if} \ x\in J_{i}.
\]
We def\/ine
\[
\Omega _{f}:=\Big\{x\in I: \ f^{k}(x)\in\bigcup _{j=1}^{m}J_{j}\ \text{for all} \ k=0,1,\dots\Big\}.
\]

Note that $\overline{\Omega }_{f}=I$. The \textit{itinerary map}
${\rm it}: \Omega _{f}\rightarrow \{1,2,\dots,m\}^{\mathbb{N}}$ is def\/ined
by
\[
{\rm it}(x)={\rm ad}(x){\rm ad}(f(x)){\rm ad}\big(f^{2}(x)\big)\cdots
\]
and let
\begin{equation}
\Lambda_{f}={\rm it}(\Omega _{f}). \label{aperiodicshift}
\end{equation}
The space $\Lambda_{f}$ is invariant under the \textit{%
shift map} $\sigma :\{1,2,\dots,m\}^{\mathbb{N}}\rightarrow
\{1,2,\dots,m\}^{\mathbb{N}}$ def\/ined by
\[
\sigma (i_{1}i_{2}\cdots)=(i_{2}i_{3}\cdots),
\]%
and we have ${\rm it}\circ f=\sigma \circ {\rm it}.$ We will use~$\sigma$
meaning in fact $\sigma _{|\Lambda_{f}}$. A sequence in
$\{1,2,\dots,m\}^{\mathbb{N}}$ is called \textit{admissible}, with
respect to~$f$, if it occurs as an itinerary for some point $x$ in~$I$, that is, if it belongs to~$\Lambda_{f}$. An
\textit{admissible\ word} is a f\/inite sub-sequence of some
admissible\ sequence. The set of admissible words of size~$k$ is
denoted by $W_{k}=W_{k}(f)$. Given $i_{1}\cdots i_{k}\in W_{k}$, we
def\/ine $I_{i_{1}\cdots i_{k}}$ as the set of points $x$ in $\Omega_f$
which satisfy
\[
{\rm ad}(x)=i_{1}, \ \dots, \ {\rm ad}\big(f^{k}(x)\big)=i_{k}.
\]%

\subsection{Linear mod 1 interval maps}

Now, let us consider the family of linear mod 1 transformations
as in equation~(\ref{linearmod1}).
In the sequel we will denote
$f_{\beta,\alpha}$ by $f$.
 The behavior of the dynamical system $(I,{f})$ is
characterized by the sequences ${\rm it}(f(c_{j}^{+}))$ and
${\rm it}(f(c_{j}^{-}))$, for each discontinuity point $c_{j}$, see~\cite{Ma-SR}. Let us consider the partition of monotonicity
$\mathcal{I}=\{I_{1},\dots,I_{n}\}$ of $f$, with
\begin{gather}
I_{1}=\left]0, {(1-\alpha)/\beta}\right[ ,  \ \dots, \ I_{j}=\left]
(j-\alpha
)/\beta ,(j+1-\alpha )/\beta \right[,\ \dots,  \nonumber \\
\hphantom{I_{1}=\left]0, {(1-\alpha)/\beta}\right[ ,}{} \ \dots, \  I_{n}=\left] (n-1-\alpha )/\beta ,1\right[,
\label{partitionm}
\end{gather}
which is the minimal partition of monotonicity for $f$
($n=[\beta]+1$ with $[\beta]$ being the integral part of $\beta$ and $\alpha>0$.
For $\alpha=0$: $n=[\beta]$ if $\beta$ is an integer, and $n=[\beta]+1$ if $\beta\notin\mathbb{N}$).
A~characterization of the values of $\alpha$, for which there is a
Markov partition, is partially given by the following:

\begin{proposition}\label{prop1} If ${\rm it}_{f}(0)=(\xi_1, \xi_2,\dots,\xi_l,\dots)$ and
${\rm it}_{f}(1)$ are periodic $($with $\xi_l=\xi_1)$ then
\[ \alpha
=\frac{\xi _{l}+\xi _{l-1}\beta+\xi _{l-2}\beta^{2}+\cdots+\xi
_{1}\beta^{l-1}}{
1+\beta+\beta^{2}+\cdots+\beta^{l-1}}.
\] In particular $\alpha \in \mathbb{Q}(\beta)$.
\end{proposition}

\begin{proof} See \cite[Proposition~2.6]{Ma-SR} for full details.
\end{proof}

\section{Subshift algebras from linear mod 1 transformations}
\label{matexel}

As in \cite{CMP}, we consider the equivalence relation
\begin{equation*}
R_{f}=\{(x,y):\ f^{n}(x)=f^{m}(y)\ \text{for some} \ n,m\in \mathbb{N}_{0}\}
.  
\end{equation*}
We write $x\sim y$ whenever $(x,y)\in R_{f}$. Consider the
equivalence class $R_{f}(x)$ ($=\bigcup_{j\in\mathbb{Z}}f^j(x)$
also called the generalized orbit of $x$) and set $H_{x}$ the
Hilbert space $l^{2}(R_{f}(x))$ with canonical orthonormal basis
$\{\left\vert y\right\rangle :y\in R_{f}(x)\}$, in Dirac notation.
Note that $H_{x}=H_{y}$ (are the same Hilbert spaces) whenever
$x\sim y$. The inner product $(\cdot ,\cdot )$ is given by
\[
\left\langle y|z\right\rangle =\left(\left\vert y\right\rangle
,\left\vert z\right\rangle \right)= \delta _{y,z}.
\]

Let now $f$ be the linear mod 1 transformation def\/ined in equation~(\ref{linearmod1}) and $\mathcal{I}=\{I_{1},\dots,I_{n}\}$ be the
partition of monotonicity as written down in equation~(\ref{partitionm}). For every $i=1,\dots,n$, let $f_{i}:=f_{|I_{i}}$
be the restriction of $f$ to the subinterval $I_i$ of the
partition $\mathcal{I}$. For every $i\in \left\{1,2,\dots,n\right\}$
let us def\/ine an operator $T_{i}$ on $H_{x}$ def\/ined f\/irst on the orthonormal basis as follows:
\begin{equation}
T_{i}\left\vert y\right\rangle =\chi _{f(I_{i})}(y)\left\vert
f_{i}^{-1}(y)\right\rangle \label{matsumotot}
\end{equation}
and then extend it by linearity and continuity to~$H_{x}$. Note
that $\chi _{f(I_{i})}(x)=1$ if and only if there is a pre-image
of~$x$ in $I_{i}$. The we have
\begin{equation}
T_{i}^{\ast }\left\vert y\right\rangle =\chi _{I_{i}}(y)\left\vert
f(y)\right\rangle. \label{adjoint}
\end{equation}
Indeed, on one hand
\[
(|y\rangle, T_{i}| z\rangle) =\big(|y\rangle,
\chi_{f(I_i)}|f_{i}^{-1}(z)\rangle\big)
\]
and on the other hand we have
\[
(T_{i}^{\ast }|y\rangle, |z\rangle)=
(\chi_{I_i}(y)|f(y)\rangle, |z\rangle)
=\chi _{I_{i}}(y)\delta
_{f(y),z}\text{.}
\]
So since $\chi_{f(I_i)}|f_{i}^{-1}(z)=\chi _{I_{i}}(y)\delta
_{f(y),z}$ we have shown that the adjoint of $T_{i}$ is given by
equation~(\ref{adjoint}). We further remark that $T_{i}$ is a partial
isometry: namely, $T_{i}$ is an isometry on its restriction to
span$\{\left\vert y\right\rangle: y\in {f(I_{i})\}\cap H}_{x}$ and
vanishes in the remaining part of $H_{x}$.

\begin{lemma}\label{CK2}
The operators $T_{i}$ satisfy the relations
\[
\sum_{i=1}^{n}T_{i}T_{i}^{\ast }={\bf 1},\qquad T_{\mu }^{\ast
}T_{\mu }T_{\nu }=T_{\nu }T_{\mu \nu }^{\ast }T_{\mu \nu }\qquad
\hbox{and}\qquad T_\mu^\ast T_\mu T_\nu^\ast T_\nu=T_\nu^\ast T_\nu
T_\mu^\ast T_\mu,
\]
 for $\mu$, $\nu $ given admissible words.
\end{lemma}
\begin{proof} Consider $T_{i}T_{i}^{\ast }$ acting on a vector
$\left\vert y\right\rangle$ of the canonical basis of $H_{x}$,
\begin{gather*}
T_{i}T_{i}^{\ast }\left\vert y\right\rangle =
\chi_{I_{i}}(y)T_{i}\left\vert f(y)\right\rangle
 = \chi _{I_{i}}(y)\chi _{f(I_{i})}(f(y))\left\vert
f_{i}^{-1}\circ
f(y)\right\rangle
 = \chi _{I_{i}}(y)\left\vert y\right\rangle,
\end{gather*}
since $\chi _{I_{i}}(y)\chi _{f(I_{i})}(f(y))=\chi _{I_{i}}(y)$.
Then
\[
(T_{1}T_{1}^{\ast }+\cdots +T_{n}T_{n}^{\ast})\left\vert
y\right\rangle =(\chi _{I_{1}}(y)+\cdots +\chi
_{I_{n}}(y))\left\vert y\right\rangle =\left\vert
y\right\rangle.
\]
 Now, consider $T_{\mu }^{\ast }T_{\mu }T_{\nu }$
acting on a vector $\left\vert y\right\rangle $ of the canonical
basis for some $\mu =\mu _{1}\cdots \mu _{k}$, $\nu =\nu
_{1}\cdots \nu_{r}$ admissible words,
\begin{gather*}
T_{\mu }^{\ast }T_{\mu }T_{\nu }\left\vert y\right\rangle
 = T_{\mu }^{\ast }T_{\mu} \chi_{f^{r}(I_{\nu })}(y)\left\vert f_{\nu
_{1}}^{-1}\circ \cdots \circ f_{\nu _{r}}^{-1}(y)\right\rangle  \\
\hphantom{T_{\mu }^{\ast }T_{\mu }T_{\nu }\left\vert y\right\rangle}{}
= T_{\mu }^{\ast }\chi _{f^{k+r}(I_{\mu \nu })}(y)\left\vert
f_{\mu _{1}}^{-1}\circ \cdots \circ f_{\mu _{k}}^{-1}\circ f_{\nu
_{1}}^{-1}\circ \cdots \circ f_{\nu _{r}}^{-1}(y)\right\rangle  \\
\hphantom{T_{\mu }^{\ast }T_{\mu }T_{\nu }\left\vert y\right\rangle}{}
=\chi _{f^{k+r}(I_{\mu \nu })}(y)\left\vert
f_{\nu_{1}}^{-1}\circ \cdots \circ f_{\nu _{r}}^{-1}(y)\right\rangle
.
\end{gather*}

On the other hand
\begin{gather}
T_{\nu }T_{\mu \nu }^{\ast }T_{\mu \nu }\left\vert y\right\rangle
= T_{\nu }T_{\mu \nu }^{\ast }\chi _{f^{k+r}(I_{\mu \nu
})}(y)\left\vert f_{\mu _{1}}^{-1}\circ \cdots \circ f_{\mu _{k}}^{-1}\circ
f_{\nu _{1}}^{-1}\circ \cdots \circ f_{\nu
_{r}}^{-1}(y)\right\rangle \nonumber\\ 
\hphantom{T_{\nu }T_{\mu \nu }^{\ast }T_{\mu \nu }\left\vert y\right\rangle}{}
 = T_{\nu }\chi _{f^{k+r}(I_{\mu \nu })}(y)\left\vert
y\right\rangle  
 = \chi _{f^{k+r}(I_{\mu \nu })}(y)\left\vert f_{\nu
_{1}}^{-1}\circ \cdots \circ f_{\nu _{r}}^{-1}(y)\right\rangle.\label{m3}
\end{gather}
Finally since $T_{\mu }^{\ast }T_{\mu }|y\rangle =\chi
_{f^{k}(I_{\mu })}(y)|y\rangle $ and $T_{\nu }^{\ast }T_{\nu
}|y\rangle =\chi _{f^{r}(I_{\nu })}(y)|y\rangle $\ for admissible
words $\mu$, $\nu $, we easily conclude that $T_{\mu }^{\ast }T_{\mu
}T_{\nu }^{\ast }T_{\nu }=T_{\nu }^{\ast }T_{\nu}T_{\mu }^{\ast
}T_{\mu }$.
\end{proof}

 As an immediate consequence of Lemma \ref{CK2}
(and above equation~(\ref{m3})) we obtain the
following.

\begin{proposition}\label{matrep}
Let $\mathcal{O}_{\Lambda_f}$ be the subshift algebra associated
to the subshift $\Lambda_f$ as defined in \eqref{aperiodicshift}
above, then $\rho_x: \mathcal{O}_{\Lambda_f}\to B(H_x)$ defined by
$t_i\to T_i$ is a representation of $\mathcal{O}_{\Lambda_f}$.
\end{proposition}

We remark here that $T_i^\ast T_i={\bf 1}$ for all $i=2, \dots , n-1$ and
$T_1^\ast T_1={\bf 1}$ if and only if $\alpha=0$. Besides $T_n^\ast
T_n={\bf 1}$ if and only if $\beta=n\in\mathbb{N}$. Therefore $\rho_x$
is a representation of a Cuntz algebra if and only if $\alpha=0$
and $\beta$ is a positive integer. In this case, the interval map
$f$ is a Markov map and moreover the partition with the
monotonicity intervals, as in (\ref{partitionm}), reduces to
\[
I_1={}]0,1/n[,\quad I_2={}]1/n,2/n[, \quad \dots ,\quad  I_{j}={}
](j-1)/n,j/n[, \quad \dots , \quad I_n={}](n-1)/n,1[
\]
and coincides with the (minimal)
Markov partition~\cite{CMP}.
From the viewpoint of interval maps, these Cuntz algebra~$\mathcal{O}_n$ representations were treated in \cite[Remark~2.9]{CMP}.


We remark that if
$f_{\beta, \alpha}$ is a linear mod~1 map with $\alpha\notin\mathbb{Q}$ and
$\beta=1$ then $f$ is not a Markov map by Proposition~\ref{prop1}
and thus the representation $\rho_x$ of Proposition~\ref{matrep} is never a
representation of a Cuntz--Krieger algebra.

\subsection{Irreducibility of the representations}

For the linear mod $1$ transformation map $f$ and the linear
operators $T_1,\dots, T_n\in B(H_x)$ def\/ined in equation~(\ref{matsumotot}), we may consider the following operator
\begin{equation}
V=T_1^\ast+\cdots +T_n^\ast, \label{v1n}
\end{equation}
which satisf\/ies $V|y\rangle=|f(y)\rangle$ on every vector basis
$|y\rangle$. In general $V$ is not unitary (unless $\beta=1$ so
that $f$ becomes an invertible function). Let $U$ be the diagonal
operator
\begin{equation}
U\left\vert y\right\rangle =e^{2\pi iy}\left\vert y\right\rangle
, \label{u1n}
\end{equation}
which is an unitary operator, with $U^{\ast }\left\vert
y\right\rangle =e^{-2\pi iy}\left\vert y\right\rangle \text{.}$ In
order to emphasize that $U, V\in B(H_x)$, we write $U_x$ and
$V_x$ for the above operators $U$ and $V$, respectively.

For a self-adjoint set of operators $\mathcal{S}\subseteq B(H)$
on some Hilbert space $H$, containing the identity ${\bf 1}$, the von
Neumann algebra generated by $\mathcal{S}$ equals the double commutant
$\mathcal{S}^{\prime\prime}$, which in turn is also equal to the
closure of $\mathcal{S}$ under the strong operator topology (this
is the famous bicommutant von Neumann theorem e.g.\ the textbook~\cite{Ped}). We note that  $\mathcal{S}^{\prime}=\{t\in B(H):
ts=st,\ \hbox{for\ all}\ s\in\mathcal{S}\}$ and
$\mathcal{S}^{\prime\prime}=(\mathcal{S}^{\prime})^\prime$  and
$\mathcal{S}^{\prime\prime\prime}=\mathcal{S}^{\prime}$. Also
$s_i$ converges to $s$ in the strong operator topology if
$||(s_i-s)\xi||\to 0$ for every vector $\xi\in H$.

Let $\mathcal{A}_{\beta,\alpha}=C^\ast(U,V)$ be the $C^*$-subalgebra
of $B(H_x)$ generated by $U$ and $V$, and consider the
representation $\rho_x$ from Proposition \ref{matrep}.

\begin{proposition}
We have $\mathcal{A}_{\beta,\alpha}\subseteq
\rho_x(\mathcal{O}_{\Lambda_f})^{\prime\prime}$.
\label{pro23}
\end{proposition}
\begin{proof}\sloppy By def\/inition of the $C^*$-algebra
$\mathcal{A}_{\beta,\alpha}$, we only need to prove that~$U$,~$V$
belong to~$\rho_x(\mathcal{O}_{\Lambda_f})^{\prime\prime}$. It is
clear that $V\in \rho_x(\mathcal{O}_{\Lambda_f})\subseteq
\rho_x(\mathcal{O}_{\Lambda_f})^{\prime\prime}$. We now show that
$U\in \rho_x(\mathcal{O}_{\Lambda_f})^{\prime\prime}$. For each
$\mu \in W_{k}$, let $m(\mu )$ be some point in $I_{\mu }\cap
R_{f}(x)$. Note that if we have ${\rm it}(y)=(\alpha _{j})_{j=1}^{\infty
}$, for some point $y\in R_{f}(x)$ then $\lim\limits_{j\rightarrow \infty
}m(\alpha_{1}\cdots \alpha _{j})=y$, and the limit is independent on
the choice of $m(\alpha _{1}\cdots \alpha_{j})\in
I_{\alpha_{1}\cdots \alpha_{j}}$, since for each $j\in \mathbb{N}$
there is $r>j$ so that $I_{\alpha _{1}\cdots \alpha_{j}}\supset
I_{\alpha _{1}\cdots \alpha _{r}}$. Let $M_{k}=\sum_{\mu \in
W_{k}}e^{2\pi im(\mu )}T_{\mu }T_{\mu }^{\ast }$. We can see that
$\lim\limits_{k\rightarrow \infty }M_{k}=U$ in the strong topology, since
$\lim\limits_{k\rightarrow \infty }\left\Vert M_{k}v-Uv\right\Vert =0$,
for every $v\in H_{x}$. Therefore, $U$ is in the von Neumann
algebra generated by the operators $T_{1},\dots ,T_n$. \end{proof}

\begin{lemma}\label{lemmadistinct}
 Let $x\in I$ and $f$ be the linear {\rm mod~1}
transformation \eqref{linearmod1} with fixed $\alpha$ and
$\beta$. Let $Q\in B(H_x)$ be an operator commuting with both $U$ and $V$, then
$Q=\lambda I$ for some $\lambda\in\mathbb{C}$.
\end{lemma}

\begin{proof}  First of all we remark that the
\begin{equation}
 e^{2\pi i z}\quad \hbox{are all distinct},
\label{disctint}
\end{equation}
with $z\in R_f(x)$, since $z\in [0,1]$.
Let $Q\in \mathcal{B}(H_x)$ commuting with both $U$ and
$V$. For each $z\in R_f(x)$ let $\mu_z:=e^{2\pi i z}$. By def\/inition
$U|z\rangle=\mu_z| z\rangle$, i.e., every $\mu_z$ is an
eigenvalue of~$U$.
We easily get
\begin{equation}
U Q| z\rangle=\mu_z  Q| z\rangle \label{eq23}
\end{equation}
by applying the def\/initions and the fact that $U$ and $Q$ commute.
For every $z\in R_f(x)$, set $\xi_z= Q| z\rangle$. Then we can
write equation~(\ref{eq23}) as follows: $U\xi_z = \mu_z\xi_z$.
Since $\{|z\rangle, z\in R_f(x)\}$ is an o.n.\ basis of
$H_x$, there are constants $c_w$ with $w\in R_f(x)$ such that $\xi_z=\sum c_w
| w\rangle$. Since
\[
U\xi_z=\sum_{w\in R_f(x)} \mu_w c_w | w\rangle
\]
and
\[
\mu_z\xi_z =\sum_{w\in R_f(x)} \mu_z c_w | w\rangle,
\]
we conclude that $c_w=0$ for all $w\not= z$, because the
$\mu_w$'s are all distinct by (\ref{disctint}) and
$\{| w\rangle\}$ is an o.n.\ basis of $H_x$. It follows that
$\xi_z=c_z| z\rangle$ or equivalently
$Q| z\rangle=c_z| z\rangle$ for some $c_z\in\mathbb{C}$.
But $V| z\rangle=|f(z)\rangle$, so $VQ=QV$ gives
$c_z=c_{f(z)}$. Therefore $Q$ is a multiple of the identity
operator.
\end{proof}

\begin{theorem}\label{nicethm11}
The representation $\rho_x$ of the
subshift $C^*$-algebra $\mathcal{O}_{\Lambda_f}$ as in Proposition~{\rm \ref{matrep}} is irreducible. Moreover $\rho_x\sim \rho_y$ if and
only if $x\sim y$.
\end{theorem}

\begin{proof} To prove that $\rho_x$ is irreducible we prove that
$\rho_x(\mathcal{O}_{\Lambda_f})^\prime=\mathbb{C}I$. First note
that from Proposition \ref{pro23} we have $C^\ast(U,V)\subseteq
\rho_x(\mathcal{O}_{\Lambda_f})^{\prime\prime}$ and so taking
commutant and using von Neumann bicommutant theorem, we conclude
that $\rho_x(\mathcal{O}_{\Lambda_f})^{\prime}\subseteq
C^\ast(U,V)^\prime$.

Now let $Q\in \rho_x(\mathcal{O}_{\Lambda_f})^\prime$. So we
can conclude that $Q\in C^\ast(U,V)^\prime$ and thus $Q$ commutes
with both $U$ and $V$. By Lemma \ref{lemmadistinct} we conclude
that $Q=\lambda I$ for some $\lambda\in\mathbb{C}$. Therefore
$\rho_x(\mathcal{O}_{\Lambda_f})^\prime=\mathbb{C}I$ and so
$\rho_x$ is an irreducible representation of
$\mathcal{O}_{\Lambda_f}$.

It is clear that if $R_f(x)=R_f(y)$, then $\rho_x$ and $\rho_y$ are
unitarily equivalent. Notice that~$U_x$ and~$U_y$ have the same
eigenvalues if and only if~$x\sim y$. Hence~$\rho_x$ and $\rho_y$
can be unitarily equivalent only when $x\sim y$. \end{proof}

\begin{remark}
If $\beta=n$ is an integer and $\alpha=0$ then the partition of $I$ into
monotonicity subintervals of $f_{n,0}$ is a Markov partition, the subshift
$\Lambda_{f_{n,0}}$ is the full shift and the underlying $C^*$-algebra is the
Cuntz algebra~$\mathcal{O}_n$. Furthermore we recover in
Theorem~\ref{nicethm11} our previous result obtained in~\cite{CMP}.
\end{remark}

\subsection*{Acknowledgment}
 First author acknowledges CIMA-UE
for f\/inancial support. The other authors were partially supported
by the Funda\c{c}\~ao para a Ci\^encia e a Tecnologia through the
Program POCI 2010/FE\-DER.

\pdfbookmark[1]{References}{ref}
\LastPageEnding

\end{document}